\pdfoutput=1
\documentclass[reqno]{amsart}

\usepackage[english]{babel}
\usepackage[utf8]{inputenc}

\usepackage{multicol}
\usepackage{tabu}

\usepackage{amsmath}
\usepackage{amssymb}
\usepackage{IEEEtrantools}
\usepackage{xcolor}
\usepackage{amstext}
\usepackage{amsthm}
\usepackage{mathtools}
\usepackage{stmaryrd}
\usepackage{mathrsfs}
\usepackage{amsfonts}

\usepackage{tikz}
\usetikzlibrary{matrix,arrows,calc,fit,cd,positioning,intersections,arrows.meta,shadings,lindenmayersystems,fadings,angles,patterns,quotes,backgrounds}
\usepackage{tikz-3dplot}
\usepackage{tikz-cd}

\numberwithin{equation}{section}

\newtheorem{thm}{Theorem}[section]

\newtheorem{lem}[thm]{Lemma}

\newtheorem{ques}[thm]{Question}

\usepackage{hyperref}
\usepackage{cleveref}

\crefname{thm}{Theorem}{Theorems}

\newcommand{\A}{\mathcal{A}}
\newcommand{\R}{\mathbb{R}}
\newcommand{\C}{\mathbb{C}}
\newcommand{\Z}{\mathbb{Z}}

\renewcommand{\S}{\mathfrak{S}}
\DeclareMathOperator{\NC}{NC}
\DeclareMathOperator{\Fix}{Fix}

\begin{document}

\markboth{\hfill{\rm Giovanni Paolini and Mario Salvetti} \hfill}{\hfill {\rm The $K(\pi, 1)$ conjecture for affine Artin groups \hfill}}

\title{The $K(\pi, 1)$ conjecture for affine Artin groups}

\author{Giovanni Paolini and Mario Salvetti}

\begin{abstract}
In this summary paper, we present the key ideas behind the recent proof of the $K(\pi, 1)$ conjecture for affine Artin groups, which states that complements of locally finite affine hyperplane arrangements with real equations and stable under orthogonal reflections are aspherical. We survey three facets of the argument: the combinatorics of noncrossing partition posets associated with Coxeter groups; the appearance of dual Artin groups and the question of their isomorphism with standard Artin groups; the topological models and their interplay in the proof.
\end{abstract}

\maketitle

\setcounter{tocdepth}{1}
\tableofcontents

\section{Introduction}

This paper is meant as a brief introduction to the proof of the $K(\pi, 1)$ conjecture for affine Artin groups, published by the authors in \cite{paolini2021proof}.
The main result can be stated without any reference to Artin groups, using minimal prerequisites from algebraic topology:

\begin{thm}[\cite{paolini2021proof}]
    Let $\A$ be a locally finite arrangement of affine hyperplanes in~$\C^n$ with real equations such that, whenever $H, H' \in \A$, also the orthogonal reflection of $H$ with respect to $H'$ is in $\A$.
    Then, the complement
    \[ \C^n \setminus \bigcup_{H \in \A} H \]
    is a $K(\pi, 1)$ space (or Eilenberg-MacLane space, or classifying space).
    \label{thm:main}
\end{thm}

The $K(\pi, 1)$ conjecture is an open problem which goes back to work by Arnol'd, Brieskorn, Pham, and Thom in the 1960s (see e.g.~\cite{brieskorn1973groupes,van1983homotopy,godelle2012basic,paris2012k}).
The conjecture emerged at the crossroads of singularity theory and the topology of complements of hyperplane arrangements, %
and can be stated for arbitrary Coxeter and Artin groups.
A \textit{Coxeter group} is a group with a presentation of the following form:
\[
  W = \langle S \mid s^2 = 1 \;\; \forall s\in S \, \text{ and} \!\! \underbrace{stst\dots}_{m(s,t) \text{ letters}} \!\!=\!\! \underbrace{tsts\dots}_{m(s,t) \text{ letters}} \!\! \forall s \neq t  \rangle,
\]
where $S$ is a finite set and $m(s, t) = m(t, s) \in \{2, 3, 4, \dots, \infty\}$ for all $s \neq t$.
For instance, if $m(s, t) = 2$, then the two generators $s$ and $t$ satisfy the relation $st = ts$, i.e., they commute.
If $m(s, t) = 3$, we obtain the famous \textit{braid relation} $sts = tst$.
If $m(s, t) = \infty$, we conventionally mean that there is no relation involving $s$ and $t$.
The same data is used to define the \textit{Artin group} associated with the Coxeter group~$W$:
\[
  G_W = \langle S \mid \!\! \underbrace{stst\dots}_{m(s,t) \text{ letters}} \!\!=\!\! \underbrace{tsts\dots}_{m(s,t) \text{ letters}} \!\! \forall s, t \in S \rangle.
\]
The cardinality of $S$ is called the \textit{rank} of $W$ (or $G_W$).
Note that the definition of $G_W$ and the rank depend not only on $W$ as an abstract group, but on the presentation of $W$ as a Coxeter group; indeed, an abstract group can be a Coxeter group in multiple non-isomorphic ways (possibly with different ranks); see e.g.~\cite{bjorner2006combinatorics}.

\begin{figure}
    \centering
    \begin{tikzpicture}[scale=1, thick, black!80]
	\draw (-2.5,0) -- (2.5,0);
	\draw[rotate=60] (-2.5,0) -- (2.5,0);
	\draw[rotate=120] (-2.5,0) -- (2.5,0);
    \end{tikzpicture}
    \qquad
    \newcommand*\rows{10}
    \begin{tikzpicture}[scale=0.9, black!80, extended line/.style={shorten >=-10cm, shorten <=-10cm}, thick]
	\clip (0.35,0.6) rectangle + (5.3,5);
	\foreach \row in {-\rows, ...,\rows} {
		\draw [extended line] ($\row*(0.5, {0.5*sqrt(3)})$) -- ($(\rows,0)+\row*(-0.5, {0.5*sqrt(3)})$);
		\draw [extended line] ($\row*(1, 0)$) -- ($(\rows/2,{\rows/2*sqrt(3)})+\row*(0.5,{-0.5*sqrt(3)})$);
		\draw [extended line] ($\row*(1, 0)$) -- ($(\rows/2,{-\rows/2*sqrt(3)})+\row*(0.5,{0.5*sqrt(3)})$);
	}
    \end{tikzpicture}
    \caption{
        Reflection arrangements in $\R^2$ associated with a spherical reflection group (left) and an affine reflection group (right).
        Any reflection group admits a Coxeter presentation by taking as generating set $S$ the reflections across the hyperplanes bounding any fixed chamber.
    }
    \label{fig:reflection-arrangements}
\end{figure}
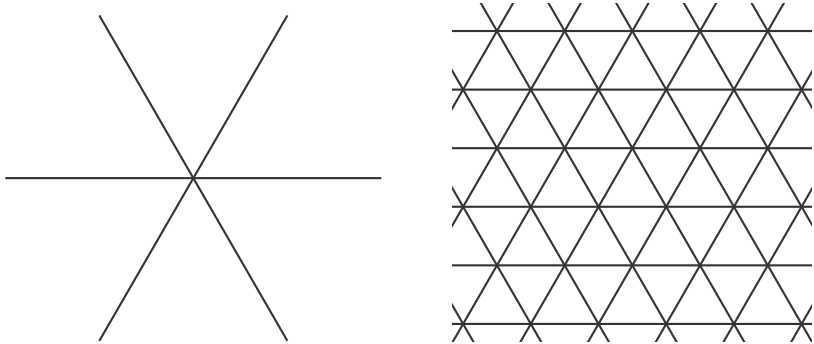

An important class of Coxeter groups is given by (real) reflection groups: discrete groups of Euclidean isometries of $\R^n$ generated by orthogonal reflections (see \Cref{fig:reflection-arrangements}).
Finite reflection groups are also called \textit{spherical}, and infinite reflection groups are also called \textit{affine}.
The corresponding Artin groups are likewise called spherical and affine, respectively.
Associated to every reflection group $W$ acting on $\R^n$ is the collection $\A$ of all hyperplanes fixed by some orthogonal reflection~$r \in W$.
On the complement $Y \subseteq \C^n$ of all complexified hyperplanes in~$\A$, the reflection group $W$ acts by a covering space action. The quotient space $Y_W = Y / W$, called the \textit{orbit configuration space} of $W$, has fundamental group isomorphic to the Artin group $G_W$; this was proved by Brieskorn in the spherical case \cite{brieskorn1971fundamentalgruppe,brieskorn1973groupes} and by van der Lek in full generality \cite{van1983homotopy}.

The previous construction can be extended to all Coxeter groups $W$, using the notion of the \textit{Tits cone} (see e.g.~\cite{bourbaki1968elements, humphreys1992reflection}): let $V$ be the standard real representation of $W$ with Tits cone $U\subseteq V$, set $U_\C = U+iV$, and define
\[ Y = U_\C \setminus \bigcup_{H\in\A} H_\C. \]
The diagonal action of $W$ on $Y$ is a covering space action.
The $K(\pi, 1)$ conjecture states that the orbit space $Y_W = Y/W$ is a $K(\pi, 1)$ space (or, equivalently, that the hyperplane complement $Y$ is).
This was proved in the spherical case by Deligne in 1972 \cite{deligne1972immeubles}, and more recently in the affine case by the authors \cite{paolini2021proof}.
\Cref{thm:main} states exactly the $K(\pi, 1)$ conjecture for (spherical and affine) reflection groups.
The ``dual approach'' introduced in \cite{paolini2021proof} was also surveyed in \cite{paolini2021dual} and extended to the rank-three case in \cite{delucchi2024dual}.

\subsection*{Acknowledgments}
It is a pleasure to thank the International Congress of Basic Science for recognizing our work on the $K(\pi,1)$ conjecture with the Frontiers of Science Award in Mathematics.
We gratefully acknowledge support from PRIN 2022A7L229 ``Algebraic and topological combinatorics'', PRIN 2022S8SSW2 ``Algebraic and geometric aspects of Lie theory'', the MIUR Excellence Department Project awarded to the Department of Mathematics of the University of Pisa (CUP I57G22000700001), and INdAM-GNSAGA.

\section{The combinatorial side: noncrossing partition posets}

Let $W$ be a Coxeter group and $S$ its generating set of reflections (usually called \textit{simple reflections}).
Denote by $R$ the set of all reflections of $W$; algebraically, $R$ can be defined as the set of all conjugates of the simple reflections.
A \textit{Coxeter element} $w \in W$ is a product of the simple reflections, each occurring exactly once, in some order.
The interval consisting of all vertices lying on shortest paths from $1$ to $w$ in the (right) Cayley graph of $W$ with respect to the generating set $R$ is called the (generalized) \textit{noncrossing partition poset} associated with the Coxeter group $W$ and the Coxeter element $w$ (as usual, it is understood that $W$ is equipped with its Coxeter group presentation).
We denote it by $\NC(W, w)$.

\begin{figure}
    \centering
    \begin{tikzpicture}[scale=1.5]
        \newcommand{\x}{1.3}
        \node (e) at (0,0) {$1$};
        \node (12) at (-1,\x) {$(12)$};
        \node (13) at (0,\x) {$(13)$};
        \node (23) at (1,\x) {$(23)$};
        \node (123) at (-1,2*\x) {$(123)$};
        \node (132) at (1,2*\x) {$(132)$};
        
        \begin{scope}[very thick]
            \draw (e.north) -- (12.south);
            \draw (e.north) -- (23.south);
            \draw (e.north) -- (13.south);
            \draw (12.north) -- (123.south);
            \draw (13.north) -- (123.south);
            \draw (23.north) -- (123.south);
        \end{scope}

        \begin{scope}[dashed]
            \draw (12.north) -- (132.south);
            \draw (13.north) -- (132.south);
            \draw (23.north) -- (132.south);            
        \end{scope}
    \end{tikzpicture}
    \qquad
    \begin{tikzpicture}[scale=1.5]
        \newcommand{\h}{1.3}
        \node (M) at (0, \h) {{$w$}};
        \node (m) at (0,-\h) {{$1$}};
        
        \foreach \x in {-5,-3,...,5} {
            \node (\x) at (\x*0.35, 0) {};
            \draw[very thick, shorten <=0.15cm] (\x.north) -- (M.south);
            \draw[very thick, shorten <=0.15cm] (\x.south) -- (m.north);
        }
        
        \foreach \x in {-5,-3,-1,...,5} {
            \pgfmathsetmacro\u{int((\x+1)/2)}; 
            \pgfmathsetmacro\v{int(\u+1)}; 
            \node  at (\x*0.35,0) {{$a_{\u}$}};
        }

        \foreach \x in {-9, -7, 7, 9} {
            \draw[dotted] (M.south) -- (\x*0.5*0.35, \h*0.5);
            \draw[dotted] (m.north) -- (\x*0.5*0.35, -\h*0.5);
        }

        \node at (-7*0.32, 0) {$\dots$};
        \node at (7*0.32, 0) {$\dots$};
    \end{tikzpicture}
    \caption{
        On the left, the Cayley graph of the symmetric group $\S_3$ (with a single edge drawn between any two connected nodes).
        The Hasse diagram of the noncrossing partition lattice $\NC(\S_3, (123))$ is highlighted.
        On the right, the poset $\NC(W, w)$ when $W$ is the affine reflection group of type $\tilde A_1$, i.e., the isometry group of the real line $\R$ generated by reflections $a_i$ with respect to integer points $i \in \Z$.
        The chosen Coxeter element $w = a_1 a_0$ is the translation $x\mapsto x+2$ and can be written as $a_{i+1} a_i$ for all $i \in \Z$.
    }
    \label{fig:ncp}
\end{figure}
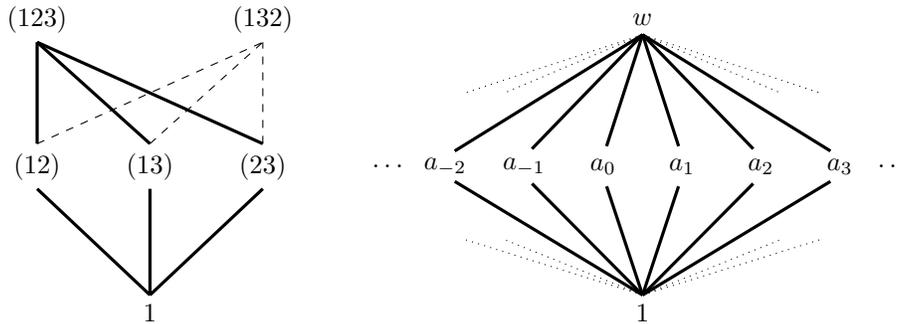

The name ``noncrossing partition poset'' originates from the case where $W$ is the symmetric group $\S_n$ and $S$ is the set of transpositions $(12), (23), \dots, (n-1 \, n)$.
In this case, Coxeter elements are $n$-cycles, and $\NC(W, w)$ is isomorphic to the classical lattice of noncrossing partitions of an $n$-gon.
This lattice has a rich and well-studied combinatorics; for instance, its elements are counted by the Catalan numbers (see e.g.\ \cite{armstrong2009generalized}).

Since $\NC(W, w)$ is an interval in the Cayley graph of $W$, the cover relations (i.e., the edges of the Hasse diagram, as in \Cref{fig:ncp}) are naturally labeled by elements of $R$.
In the spherical case, every reflection $r \in R$ appears as a label, regardless of the choice of Coxeter element $w$; in general, only the reflections in a subset $R_0 \subseteq R$ appear as labels.

\subsection{Lattice property}
\label{sec:lattice}

Spherical noncrossing partition posets all satisfy the lattice property: every pair of elements has a unique minimal upper bound and a unique maximal lower bound.
This was classically known in the case of the symmetric group, then proved for all spherical cases by Bessis with a case-by-case argument \cite{bessis2003dual} and by Brady and Watt in a uniform way \cite{brady2008non}.

Digne showed that the lattice property holds for two families of affine cases: $\tilde A_n$ (for certain choices of the Coxeter element) and $\tilde C_n$ \cite{digne2006presentations, digne2012garside}.
Together with the rank-three case $\tilde G_2$, this list turned out to be complete, as McCammond showed that the lattice property fails in all remaining cases \cite{mccammond2015dual}.

More recently, Delucchi and the authors noted that the lattice property holds for all rank-three cases \cite{delucchi2024dual}.
In rank $\geq 4$, a characterization of when $\NC(W, w)$ is a lattice is not currently known.

\subsection{Lexicographic shellability}

Another combinatorial property relevant to our work is edge-lexicographic shellability (or EL-shellability).
A bounded poset $(P, \leq\nobreak)$ (i.e., a poset with a minimum and a maximum element) is \textit{EL-shellable} if its cover relations can be labeled by elements of a totally ordered set $(\Lambda, \preceq)$ in such a way that every interval $[x, y] \subseteq P$ contains a unique inclusion-maximal chain $x = x_0 < x_1 < \dots < x_k = y$ with strictly $\preceq$-increasing labels, and this chain $\preceq$-lexicographically precedes all other inclusion-maximal chains in $[x, y]$.
Introduced by Bj\"orner in the 1980s \cite{bjorner1980shellable, bjorner1983lexicographically}, EL-shellability is a well-studied property in combinatorial topology, with strong topological consequences for the order complex of a poset (see e.g.\ \cite{wachs2006poset,kozlov2007combinatorial}).

Noncrossing partition posets come with a natural edge labeling with label set $\Lambda = R_0 \subseteq R$.
Then one can ask the following question:
\begin{ques}
    Is there a total order $\preceq$ on the set of reflections $R_0$ which makes $\NC(W, w)$ EL-shellable?
    \label{q:shellability}
\end{ques}

In the case of the symmetric group $\S_n$ with Coxeter element $w = (12\, \cdots\, n)$, this can be achieved by ordering $R_0$ (the transpositions) lexicographically:
given $i< j$ and $i' < j'$, define $(ij) \prec (i'j')$ if and only if $i < i'$ or $i = i'$ and $j < j'$.
The unique maximal chain in an interval $[x, y] \subseteq \NC(\S_n, w)$ can be found inductively, starting with $x < x \cdot (ij)$ where $i$ is the smallest non-fixed point of $x^{-1}y$ and $j = x^{-1}y(i)$.
This construction was generalized to all spherical Coxeter groups by Atha\-na\-sia\-dis, Brady, and Watt \cite{athanasiadis2007shellability}.

\begin{figure}
    \centering
    \begin{tikzpicture}[scale=1.2]
        \clip (-1.2, -1.5) rectangle (1.2, 3.5);
        
        \draw (0, -3) -- (0, 5);
        \draw[very thick, black!80] (0, 0) -- (0, 1);

        \node (basepoint) at (-0.1, 0.5) {};
        \draw[purple, -latex, very thick] (basepoint) -- ($(basepoint) + (0, 1.1)$);

        \foreach \i in {-5, ..., 5} {
            \node[draw, circle, fill=black, inner sep=1, label={[xshift=0.1pt]right:$a_{\i}$}] at (0, \i) {};
        }

    \end{tikzpicture}
    \qquad
    \begin{tikzpicture}[scale=1.2,
        extended line/.style={shorten >=-#1,shorten <=-#1}, extended line/.default=35cm]
        
        \clip (-2.2, -2.7) rectangle (3.1, 2.2);
        
        \begin{scope}[every path/.style={black!25}]
            \fill (0, -10) -- (0, 10) -- ($({0.5*sqrt(3)}, 10)$) -- ($({0.5*sqrt(3)}, -10)$) -- cycle;
        \end{scope}
        
        \fill[black!50] ($(0, 0)$) -- ($(0, -1)$) -- ($({0.5*sqrt(3)}, -0.5)$) -- cycle;

        \draw [extended line, dashed, thick] ($({0.25*sqrt(3)}, 0)$) -- ($({0.25*sqrt(3)}, 1)$);

        \begin{scope}
            \newcommand{\rows}{5}
            \foreach \row in {-\rows, ...,\rows} {
                \draw [extended line] ($\row*({0.5*sqrt(3)}, 0.5)$) -- ($(0,\rows)+\row*({0.5*sqrt(3)}, -0.5)$);
                \draw [extended line] ($\row*(0, 1)$) -- ($({\rows/2*sqrt(3)}, \rows/2)+\row*({-0.5*sqrt(3)}, 0.5)$);
                \draw [extended line] ($\row*(0, 1)$) -- ($({-\rows/2*sqrt(3)}, \rows/2)+\row*({0.5*sqrt(3)}, 0.5)$);
            }
        \end{scope}
        
        \node (basepoint) at ($({0.25*sqrt(3)}, -0.5)$) {};
        \draw[purple, -latex, very thick] (basepoint) -- ($(basepoint) + (0, 0.6)$);

        \node (a) at (0.65,-0.18) {$a_1$};
        \node (c) at (0.65,-0.84) {$c_0$};
        \node (c2) at (0.65,0.16) {$c_2$};
        \node (a3) at (0.65,0.82) {$a_3$};
        \node (a5) at (0.65,1.82) {$a_5$};
        \node (c4) at (0.65,1.16) {$c_4$};

        \coordinate (a-1-pos) at (0.72,-1.18);
        \coordinate (c-2-pos) at (0.72,-1.84);
        \coordinate (a-3-pos) at (0.72,-2.18);

        \begin{scope}
            \clip (0, -2.5) rectangle ($({0.5*sqrt(3)+0.01}, 10)$);
            \node[fill=black!25, inner sep=0, text opacity=0] at (a-1-pos) {$a_{-1}$};
            \node[fill=black!25, inner sep=0, text opacity=0] at (c-2-pos) {$c_{-2}$};
            \node[fill=black!25, inner sep=0, text opacity=0] at (a-3-pos) {$a_{-3}$};
        \end{scope}

        \node (a-1) at (a-1-pos) {$a_{-1}$};
        \node (c-2) at (c-2-pos) {$c_{-2}$};
        \node (a-3) at (a-3-pos) {$a_{-3}$};
        
        \node (b) at (-0.15,-0.5) {$b$};
        \node (bprime) at (1.05,0) {$b'$};
        
    \end{tikzpicture}
    \caption{
        Construction of the axial order in the cases $\tilde A_1$ (left) and $\tilde A_2$ (right).
        On the left, the order is $a_1 \prec a_2 \prec a_3 \prec \dots \prec a_{-1} \prec a_{0}$.
        On the right, the axis $\ell$ is dashed, and the lines corresponding to reflections in $R_0$ are those that intersect the shaded strip.
        Their order is $a_1 \prec c_2 \prec a_3 \prec \dots \prec b \prec b' \prec \dots \prec a_{-1} \prec c_0$;
        here, the reflections $b$ and $b'$ can be swapped.
    }
    \label{fig:axial-order}
\end{figure}
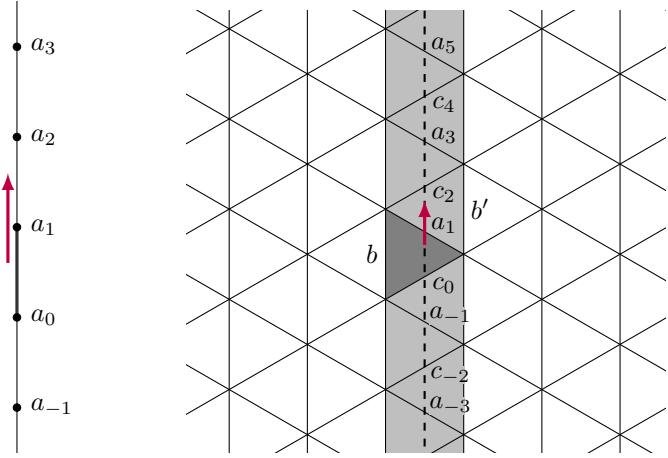

Question \ref{q:shellability} arose as a crucial step in the proof of the $K(\pi, 1)$ conjecture in the affine case.
In fact, a substantial part of the work in \cite{paolini2021proof} was devoted to constructing a suitable order of the reflections and proving EL-shellability.
The idea of the construction is quite geometric: 
an affine Coxeter element $w$ is a Euclidean isometry of $\R^n$ of hyperbolic type (i.e.,\ with no fixed points) having a $1$-dimensional line $\ell \subseteq \R^n$ of minimally moved points.
Traverse the axis $\ell$ from a fixed generic point $p_0 \in \ell$ in the direction of motion dictated by $w$, pass through infinity, and then return to $p_0$ from the opposite direction.
The total order~$\preceq$ on $R_0$ is obtained by keeping track of the order in which reflection hyperplanes are crossed.
See \Cref{fig:axial-order}.

In rank $\geq 4$, an important subtlety arises: many reflections in $R_0$ have fixed hyperplanes parallel to the axis $\ell$ (these are called \textit{horizontal reflections}), and their relative order cannot be chosen arbitrarily.
Ties are broken by infinitesimally tilting the axis $\ell$ in a suitable direction to order the horizontal reflections in a maximal finite subgroup, and then extending the order to all remaining horizontal reflections.

We call the resulting total order on $R_0$ an \textit{axial order}.
The fact that an axial order satisfies the EL-shellability property is far from obvious.
The most crucial part of the proof consists of showing that any hyperbolic element $u \in \NC(W, w)$ is a Coxeter element of a lower-rank affine Coxeter group $W_u$ and the axis of $u$ is ``sufficiently close'' to $\ell$, so that the restriction of the axial order to $R_0 \cap [1, u]$ is an axial order for $[1, u] = \NC(W_u, u)$.
For this, a key step is \cite[Lemma 3.21]{paolini2021proof}, which we quote here in its entirety:

\begin{lem}[{\cite[Lemma 3.21]{paolini2021proof}}]
    Let $W$ be an irreducible affine Coxeter group, $w$ one of its Coxeter elements, and $u \in \NC(W, w)$ a hyperbolic element such that $W_u$ is irreducible.
	Let $p$ be a point of $\ell$ that does not lie on any hyperplane of $\A_u$, and let $C$ be the chamber of $\A_u$ containing $p$.
	Then $C$ has exactly $l(u)$ walls, and $u$ can be written as the product of the reflections with respect to the walls of $C$ in the following order:
	\begin{itemize}
		\item first there are the vertical reflections that fix a point of $\ell$ above $p$, and $r$ comes before $r'$ if $\Fix(r) \cap \ell$ is below $\Fix(r') \cap \ell$;
		\item then there are the horizontal reflections, in some order;
		\item finally there are the vertical reflections that fix a point of $\ell$ below $p$, and again $r$ comes before $r'$ if $\Fix(r) \cap \ell$ is below $\Fix(r') \cap \ell$.
	\end{itemize}
\end{lem}

This lemma was checked case by case for the infinite families of affine Coxeter groups, and by computer for the exceptional groups.
A case-free proof would likely provide new insights into the fundamental reasons why the axial order yields EL-shellability of $\NC(W, w)$.

In \cite{delucchi2024dual}, Delucchi and the authors showed that the axial order provides a positive answer to Question \ref{q:shellability} for all rank-three Coxeter groups (only a finite number of which are spherical or affine).
Here, the definition of axial order was adjusted to capture the intersections ``beyond infinity'' between the axis $\ell$ and the hyperplanes fixed by horizontal reflections.
However, in work that is still in preparation, it is shown that the axial order does not always satisfy the EL-shellability property in rank~$\geq 4$.
New ideas are therefore needed to tackle Question \ref{q:shellability} in the general~case.

\section{The algebraic side: dual Artin groups}

The connection between noncrossing partition posets and Artin groups was first discovered in the case $W = \S_n$ (the corresponding Artin group $G_{\S_n}$ is the classical \textit{braid group}): Birman, Ko, and Lee \cite{birman1998new} gave a new presentation for the braid group arising from the combinatorics of noncrossing partitions.
Later, Bessis generalized this construction to all spherical Artin groups \cite{bessis2003dual}.

The resulting presentations are known as \textit{dual presentations} and are defined as follows:
the set of reflections $R_0$ serves as the abstract set of generators;
for every pair of inclusion-maximal chains $1 < r_1 < r_1r_2 < \dots <  r_1r_2\cdots r_k = w$ and $1 < r_1' < r_1' r_2' < \dots < r_1' r_2' \cdots r_k' = w$ with $r_i \in R_0$, add the relation $r_1 r_2 \cdots r_k = r_1' r_2' \cdots r_k'$.
For example, for $W = \S_3$ and $w = (123)$ as in \Cref{fig:ncp} (left), one obtains the dual presentation for the braid group on three strands:
\[ \langle a, b, c \mid ab = bc = ca \rangle, \]
where $a, b, c$ are the abstract generators corresponding to the transpositions $(12)$, $(23)$, and $(13)$.

The dual presentation can be introduced for every noncrossing partition poset $\NC(W, w)$, even in non-spherical cases.
The resulting group $G^*_{W, w}$ is called a \textit{dual Artin group}.
A key open question in the theory of Artin groups is the following:

\begin{ques}
    Is every dual Artin group $G^*_{W, w}$ isomorphic to the corresponding standard Artin group $G_W$?
    \label{q:standard-dual-isomorphism}
\end{ques}

The aforementioned result by Bessis \cite{bessis2003dual} provides a positive answer in all spherical cases.
McCammond and Sulway were able to give a positive answer in all affine cases \cite{mccammond2017artin}; as we explain in \Cref{sec:topology}, a byproduct of the proof of the $K(\pi, 1)$ conjecture is a different, topologically flavored, proof of the isomorphism $G^*_{W, w} \cong G_W$ in all affine cases.
More recently, Delucchi and the authors gave a positive answer to Question \ref{q:standard-dual-isomorphism} in the rank-three case \cite{delucchi2024dual}.

\section{The topological side: CW models and homotopy equivalences}
\label{sec:topology}

\subsection{Interval complex}

\begin{figure}
    \centering
    \begin{tikzpicture}[scale=1.5]
        \newcommand{\h}{1.3}
        \coordinate (M) at (0, \h) {};
        \coordinate (m) at (0,-\h) {};

        \foreach \x in {-3, -1, 3} {
            \coordinate (\x) at (\x*0.35, 0) {};
        }

        \fill[fill=black!30!white](-1) -- (M) -- (m) -- cycle;
        \fill[fill=black!10!white](-3) -- (M) -- (-1) -- (m.center) -- cycle;
        \fill[fill=black!20!white](3) -- (M) -- (m) -- cycle;

        \draw (3) -- node[pos=0.4, right, yshift=1pt] {\footnotesize $c$} (M);
        \draw (3) -- node[pos=0.4, right, yshift=-1pt] {\footnotesize $b$} (m);

        \draw (-3) -- node[pos=0.4, left, yshift=1pt] {\footnotesize $b$} (M);
        \draw (-3) -- node[pos=0.4, left, yshift=-1pt] {\footnotesize $a$} (m);

        \draw (-1) -- node[pos=0.4, left, xshift=1pt, yshift=1pt] {\footnotesize $a$} (M);
        \draw (-1) -- node[pos=0.4, left, xshift=1pt, yshift=-1pt] {\footnotesize $c$} (m);

        \draw (m) -- node[pos=0.5, right] {\footnotesize $w$} (M);
    \end{tikzpicture}
    \qquad\qquad
    \begin{tikzpicture}[scale=1.5]
        \newcommand{\h}{1.3}
        \coordinate (M) at (0, \h) {};
        \coordinate (m) at (0,-\h) {};

        \foreach \x in {-5,-3,...,5} {
            \coordinate (\x) at (\x*0.35, 0) {};
        }

        \fill[fill=black!30!white](-1) -- (M) -- (m) -- cycle;
        \fill[fill=black!40!white](1) -- (M) -- (m) -- cycle;
        \fill[fill=black!20!white](-3) -- (M) -- (-1) -- (m.center) -- cycle;
        \fill[fill=black!25!white](3) -- (M) -- (1) -- (m) -- cycle;
        \fill[fill=black!10!white](-5) -- (M) -- (-3) -- (m) -- cycle;
        \fill[fill=black!10!white](5) -- (M) -- (3) -- (m) -- cycle;

        \foreach \x in {1,3,5} {
            \pgfmathsetmacro\u{int((\x+1)/2)}; 
            \pgfmathsetmacro\v{int(\u-1)}; 

            \draw (\x) -- node[pos=0.25, right, xshift=0pt, yshift=1pt] {\footnotesize $a_{\v}$} (M);
            
            \draw (\x) -- node[pos=0.25, right, xshift=0pt, yshift=-1pt] {\footnotesize $a_{\u}$} (m);
        }

        \foreach \x in {-5,-3,-1} {
            \pgfmathsetmacro\u{int((\x+1)/2)}; 
            \pgfmathsetmacro\v{int(\u-1)}; 

            \draw (\x) -- node[pos=0.25, left, xshift=2pt, yshift=1pt] {\footnotesize $a_{\v}$} (M);
            
            \draw (\x) -- node[pos=0.25, left, xshift=2pt, yshift=-1pt] {\footnotesize $a_{\u}$} (m);
        }

        \draw (m) -- node[pos=0.5, left, xshift=2pt] {\footnotesize $w$} (M);

        \foreach \x in {-9, -7, 7, 9} {
            \draw[dotted] (M) -- (\x*0.5*0.35, \h*0.5);
            \draw[dotted] (m) -- (\x*0.5*0.35, -\h*0.5);
        }

        \node at (-7*0.32, 0) {$\dots$};
        \node at (7*0.32, 0) {$\dots$};
    \end{tikzpicture}
    \caption{
        Interval complex $K_{W, w}$ in the spherical case $\S_3$ (left) and in the affine case $\tilde A_1$ (right).
        They correspond to the noncrossing partition posets of \Cref{fig:ncp}.
        In either complex, all vertices are identified, and edges with the same label are identified (as oriented edges going from bottom to top).
    }
    \label{fig:interval-complex}
\end{figure}
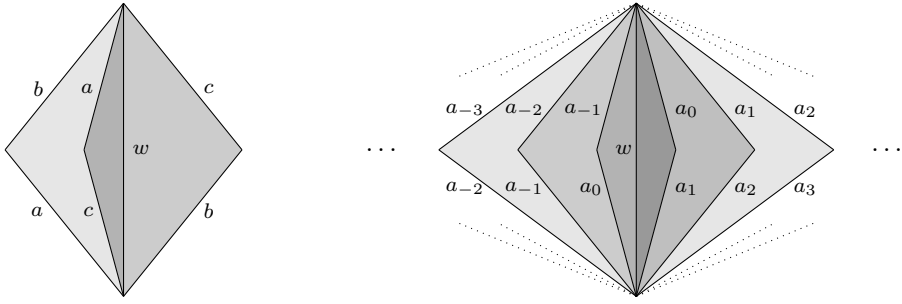

When a noncrossing partition poset $\NC(W, w)$ is a lattice, the powerful machinery of Garside groups (see e.g.~\cite{dehornoy1999gaussian, dehornoy2003homology, charney2004bestvina, dehornoy2015foundations}) allows us to explicitly construct a $\Delta$-complex $K_{W, w}$ of type $K(G^*_{W, w}, 1)$:
for every $u \in \NC(W, w)$ and for every factorization $u = x_1 x_2 \cdots x_d$ with $l(x_1) + l(x_2) + \dots l(x_d) = l(u)$ and $l(x_i) \geq 1$ (here $l(\cdot)$ denotes reflection length with respect to $R$), there is a $d$-simplex denoted by $[x_1|x_2|\dots|x_d]$; its codimension-one faces are given by the simplices $[x_2|\dots|x_d]$, $[x_1|\dots |x_i x_{i+1}| \dots | x_d]$ for $i=1, \dots, d-1$, and $[x_1|\dots | x_{d-1}]$.
We call $K_{W, w}$ the \textit{interval complex} associated with $\NC(W, w)$; see \Cref{fig:interval-complex}.
Note that the fundamental group of $K_{W, w}$ is precisely the dual Artin group $G^*_{W, w}$, by construction.

As mentioned in \Cref{sec:lattice}, the lattice property is known to hold for all spherical and rank-three Coxeter groups, as well as for some of the affine Coxeter groups.
In the proof of the $K(\pi, 1)$ conjecture for affine Artin groups, a crucial step consisted in showing that $K_{W, w}$ is in fact a $K(\pi, 1)$ for all affine noncrossing partition posets, even when the lattice property fails \cite[Theorem 6.6]{paolini2021proof}.
This opened up the following question:

\begin{ques}
    Is the interval complex $K_{W, w}$ a $K(\pi, 1)$ for all noncrossing partition posets $\NC(W, w)$?
    \label{q:interval-complex-kpi1}
\end{ques}

\subsection{Models for the orbit configuration space}

The orbit configuration space $Y_W$ is classically known to have the homotopy type of a CW complex $X_W$ with cells $\sigma_T$ indexed by the subsets $T \subseteq S$ generating a spherical subgroup $W_T$ of $W$ \cite{salvetti1987topology, salvetti1994homotopy}; this complex is known as the Salvetti complex of $W$.

Substituting each cell $\sigma_T$ of $X_W$ with a copy of the interval complex $K_{W_T, w_T}$ (for suitable Coxeter elements $w_T$ of $W_T$), the authors constructed a new CW model $X_{W, w}' \simeq Y_W$ which is a subcomplex of $K_{W, w}$.
More precisely, the definition is as follows:
\[
  X_{W, w}' \;=\; 
    \bigcup_{\substack{T\subseteq S\\[0.05cm] W_T\text{ spherical}}}
      K_{W_T,w_T}
    \;\subseteq\;
    K_{W,w}.
\]

It is worth pointing out that the construction of $X_{W, w}'$, as well as the homotopy equivalence $X_{W, w}' \simeq Y_W$, holds for all Coxeter groups \cite[Theorem 5.5]{paolini2021proof}.
The following question naturally arises:

\begin{ques}
    Is $X_{W, w}'$ homotopy equivalent to (or even a deformation retract of) the interval complex $K_{W, w}$?
    \label{q:deformation-retraction}
\end{ques}

Note that Questions \ref{q:interval-complex-kpi1} and \ref{q:deformation-retraction} together imply the $K(\pi, 1)$ conjecture.
Additionally, Question \ref{q:deformation-retraction} implies Question \ref{q:standard-dual-isomorphism}, since the complexes $X_{W, w}'$ and $K_{W, w}$ have fundamental groups isomorphic to $G_W$ and $G^*_{W, w}$, respectively.
See also \cite[Figure 8]{paolini2021dual} for a graph of implications including the previous ones.

\subsection{Deformation retraction}

A large part of \cite{paolini2021proof} is dedicated to constructing a deformation retraction $K_{W, w} \searrow X_{W, w}'$ for all affine noncrossing partition posets $\NC(W, w)$, thus answering Question \ref{q:deformation-retraction} and concluding the proof of the $K(\pi, 1)$ conjecture in the affine case.
Using the language of discrete Morse theory, the problem reduces to finding an \textit{acyclic matching} (or \textit{discrete Morse vector field}) on the face poset of $K_{W,w}$ that is \emph{complete} on $K_{W,w}\setminus X'_{W,w}$: every cell of $K_{W,w}\setminus X'_{W,w}$ is matched, so the only critical cells lie in $X'_{W,w}$.
The construction relies on a deep understanding of the combinatorics of $\NC(W, w)$, with a key role played by EL-shellability (Question~\ref{q:shellability}).

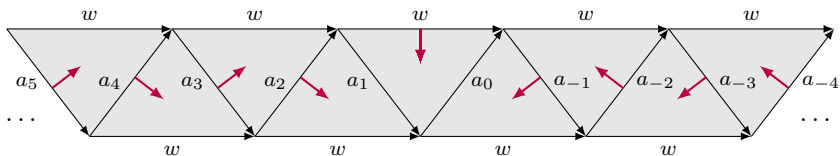
\begin{figure}
    \centering
    \begin{tikzpicture}[scale=1.1]
        \newcommand{\h}{1.3}
        \newcommand{\len}{0.35}
    
        \fill[black!10] (-4, \h) -- (6, \h) -- (5, 0) -- (-3, 0) -- cycle;
    
        \foreach \x in {-4, -2, ..., 4} {
            \pgfmathsetmacro\u{int(-\x)}; 
            \pgfmathsetmacro\v{int(\u+1)};
            \pgfmathsetmacro{\labside}{ifthenelse(\x>0,"right","left")}
            \pgfmathsetmacro{\labsideother}{ifthenelse(\x>=0,"right","left")}
            
            \draw[-latex] (\x, \h) -- node[pos=0.5, above] {\footnotesize $w$} (\x+2, \h);
            \ifnum\x<4
                \draw[-latex] (\x+1, 0) -- node[pos=0.5, below] {\footnotesize $w$} (\x+3, 0);
            \fi
            
            \draw[-latex] (\x, \h) -- node[pos=0.5, \labside] {\footnotesize $a_{\v}$} (\x+1, 0);
            \draw[-latex] (\x+1, 0) -- node[pos=0.5, \labsideother] {\footnotesize $a_{\u}$} (\x+2, \h);

            \ifnum\x<0
                \draw[-latex, purple, thick] (\x+0.55, \h*0.45) -- (\x+0.55 + \len, \h*0.45 + \len/\h);
                \draw[-latex, purple, thick] (\x+1.55, \h*0.55) -- (\x+1.55 + \len, \h*0.55 - \len/\h);
            \fi

            \ifnum\x>0
                \draw[-latex, purple, thick] (\x+0.45, \h*0.55) -- (\x+0.45 - \len, \h*0.55 - \len/\h);
                \draw[-latex, purple, thick] (\x+1.45, \h*0.45) -- (\x+1.45 - \len, \h*0.45 + \len/\h);
            \fi

            \draw[-latex, purple, thick] (1, \h) -- (1, \h-0.42);
        }

        \node at (-3.8, 0.2) {\dots};
        \node at (5.8, 0.2) {\dots};
    \end{tikzpicture}
    \caption{
        Deformation retraction $K_{W, w} \searrow X_{W, w}'$ in the affine case $\tilde A_1$.
        The interval complex $K_{W, w}$ is the same as the one shown in \Cref{fig:interval-complex} on the right, with a different visual arrangement of the cells.
        The purple arrows describe the discrete Morse vector field.
    }
    \label{fig:deformation-retraction}
\end{figure}

\Cref{fig:deformation-retraction} shows the construction in the simplest affine case, namely $\tilde A_1$.
In this case, the interval complex has one $0$-cell, infinitely many $1$-cells, namely $[w]$ and $[a_i]$ for $i \in \Z$, and infinitely many $2$-cells $[a_{i+1}|a_i]$ for $i \in \Z$.
The deformation retraction pushes the $1$-cell $[a_i]$ towards the $2$-cell $[a_i|a_{i-1}]$ if $i \geq 2$ and towards the $2$-cell $[a_{i+1}|a_i]$ if $i < 0$; in addition, the $1$-cell $[w]$ is pushed towards the $2$-cell $[a_1|a_0]$.
This leaves the $0$-cell and the $1$-cells $[a_0]$ and $[a_1]$, which constitute the subcomplex $X_{W, w}' \cong S^1 \vee S^1$.
The role of EL-shellability appears in the choice of matching $[w]$ with $[a_1|a_0]$, since $a_1a_0$ is the only $\prec$-increasing factorization (with the axial order $\prec$ visualized in \Cref{fig:axial-order}, left).

In \cite{delucchi2024dual}, a deformation retraction $K_{W, w} \searrow X_{W, w}'$ was constructed for the rank-three cases.
This recent development also relied on understanding the combinatorics of $\NC(W, w)$ and suggests that the ``dual approach'' to the $K(\pi, 1)$ conjecture could be viable for more general Coxeter groups.

\bibliographystyle{fsastyle}
\bibliography{bibliography}

\end{document}